\documentclass[12pt,leqno]{article}
\usepackage{amsmath, amssymb, indentfirst}

\newcommand{\sect}[1]{\setcounter{equation}{0}\section{#1}}

\def\epsilon{\varepsilon}

\begin{document}

\Large \noindent 
{\bf In Markov process, an extremal reversible} 

\noindent 
{\bf measure is an extremal invariant measure.}


\normalsize

\vspace*{0.8em}

\noindent 
Hiroki Yagisita 

\noindent 
(Department of Mathematics, Kyoto Sangyo University) 

\vspace*{3.2em}

\noindent 
Abstract: 

\noindent 
We consider a discrete-time temporally-homogeneous conservative Markov process.  
We show that extremality of reversible measure implies extremality of invariant measure. 
Using analogue of Dirichlet form, we modify a proof 
that in stochastic Ising model (Glauber dynamics), 
an extreme Gibbs state is an extreme invariant measure. 

\vspace*{0.8em}

\noindent 
Keyword: 

\noindent 
Choquet simplex, ergodic decomposition, detailed balance condition, mutual singularity.

\vfill

\noindent 

\newpage

\sect{Introduction}
Let $(X,\cal B)$ be a measurable space. 
Let $\{P^x\}_{x\in X}$ be a discrete-time temporally-homogeneous Markov process 
whose state-space is $(X,\cal B)$.
Let $\{p(x,E)\}_{x\in X,E\in\cal B}$ be the transition function of $\{P^x\}_{x\in X}$. 
That is, $$p(x,E) \, := \ P^x(\eta(1)\in E).$$
$\{P^x\}_{x\in X}$ is conservative, 
if and only if for any $x\in X$, $p(x,X)=1$ holds.

{\bf Definition 1 (extreme element):} \ 
Let $C$ be a subset of the set of all probability measures on $(X,\cal B)$. 
Let $m\in C$. 
Then, $m$ is said to be an {\it extreme element} of $C$, 
if and only if $[ \ t\in(0,1), \, m_0, m_1\in C, \, m=(1-t)m_0+tm_1 \ ]$ 
implies $m_0=m_1$. 
\hfill $\Box$

{\bf Definition 2 ($T$):} \ 
Let define the linear map $T$ 
from the set of all $\mathbb R$-valued bounded $\cal B$-measurable functions 
to the set of all $\mathbb R$-valued bounded $\cal B$-measurable functions 
by $$(Tf)(x) \, := \ \int_{y\in X}f(y)p(x,dy).$$
\hfill $\Box$

{\bf Definition 3 (${\cal I, R, G, I}_e, {\cal R}_e, {\cal G}_e$):} \ 
(1) \ Let $m$ be a probability measure on $(X,\cal B)$. 
Then, $m$ is said to be an {\it invariant measure}, if and only if 
for any $\mathbb R$-valued bounded $\cal B$-measurable function $f$, 
$$\int_{x\in X}(Tf)(x)m(dx)=\int_{x\in X}f(x)m(dx)$$
holds. $m$ is said to be a {\it reversible measure}, if and only if 
for any $\mathbb R$-valued bounded $\cal B$-measurable functions $f$ and $g$, 
$$\int_{x\in X}(Tf)(x)g(x)m(dx)=\int_{x\in X}f(x)(Tg)(x)m(dx)$$
holds. $m$ is said to be a {\it conservative measure}, if and only if 
for any $\mathbb R$-valued bounded $\cal B$-measurable function $g$, 
$$\int_{x\in X}(T1)(x)g(x)m(dx)=\int_{x\in X}g(x)m(dx)$$
holds. 

(2) \ Let define ${\cal I}$ as the set of {\it all invariant measures}. 
Let define ${\cal I}_e$ as the set of all extreme elements of ${\cal I}$. 
Let define ${\cal R}$ as the set of {\it all reversible measures}. 
Let define ${\cal R}_e$ as the set of all extreme elements of ${\cal R}$. 
Let define ${\cal G}$ as the set of {\it all conservative reversible measures}. 
Let define ${\cal G}_e$ as the set of all extreme elements of ${\cal G}$. 
\hfill $\Box$

{\bf Theorem 4 (Main Result):} \ 
{\it Suppose that $\{P^x\}_{x\in X}$ is conservative. 
Then, ${\cal R}_e\subset{\cal I}_e$ holds.}  
\hfill $\Box$

\newpage

\sect{Proof}
{\bf Definition 5 (joint distribution $\sigma_m$):} \ 
For each probability measure $m$ on $(X,\cal B)$, 
let define $\sigma_m$ as the joint distribution with $m$ as the initial distribution 
and $\{0,1\}$ as the set of times. That is, 
$$\sigma_m(C) \, := \ \int_{x_0\in X}\left(\int_{x_1\in X}1_C(x_0,x_1)p(x_0,dx_1)\right)m(dx_0).$$
\hfill $\Box$

{\bf Lemma 6:} \ 
{\it Let $m\in\cal R$. Then, $\sigma_m$ is symmetric.} 

{\it Proof:} \ Because of 
$\sigma_m(A\times B)=\int_{x_0\in X}(\int_{x_1\in X}1_A(x_0)1_B(x_1)p(x_0,dx_1))m(dx_0)$ 

\noindent 
$=\int_{x_0\in X}1_A(x_0)(T1_B)(x_0)m(dx_0)
=\int_{x_0\in X}1_B(x_0)(T1_A)(x_0)m(dx_0)
=\int_{x_0\in X}$

\noindent 
$(\int_{x_1\in X}1_B(x_0)1_A(x_1)p(x_0,dx_1))m(dx_0)
=\sigma_m(B\times A)$, it is symmetric (by Hopf extension theorem). 
\hfill $\blacksquare$

The following lemma is the basis in this section. 
It was inspired by the formula (1.4.8) in [1]. 

{\bf Lemma 7 (quadratic form):} \ 
{\it Let $m\in\cal G$. Then, 
for any $\mathbb R$-valued bounded $\cal B$-measurable function $f$, 
$$\int_{x\in X} ((1-T)f)(x)f(x)m(dx)$$
$$= \ \frac{1}{2}\int_{(x_0,x_1)\in X^2}|f(x_0)-f(x_1)|^2\sigma_m(d(x_0,x_1))$$
holds.} 

{\it Proof:} \ From Lemma 6, 
$2\int_{x\in X} ((1-T)f)(x)f(x)m(dx)
=2\int_{x\in X} f(x)$

\noindent 
$(f(x)(T1)(x)-(Tf)(x))m(dx)
=2\int_{x_0\in X} f(x_0) ( \int_{x_1\in X} (f(x_0)-f(x_1)) p(x_0,dx_1) )$ 

\noindent 
$m(dx_0)
=2\int_{(x_0,x_1)\in X^2} f(x_0) (f(x_0)-f(x_1)) \sigma_m(d(x_0,x_1))
=\int_{(x_0,x_1)\in X^2} f(x_0)$ 

\noindent 
$(f(x_0)-f(x_1))\sigma_m(d(x_0,x_1))
+\int_{(x_0,x_1)\in X^2} f(x_1) (f(x_1)-f(x_0)) \sigma_m(d(x_0,x_1))
=\int_{(x_0,x_1)\in X^2}|f(x_0)-f(x_1)|^2\sigma_m(d(x_0,x_1))$
holds. 
\hfill $\blacksquare$

{\bf Lemma 8:} \ 
{\it Let $m\in\cal G$. Let $f$ be a $\mathbb R$-valued bounded $\cal B$-measurable function. 
Then, the followings {\rm (a), (b), (c)} and {\rm (d)} are equivalent. 

{\rm (a)} \ Let $g$ be a $\mathbb R$-valued bounded $\cal B$-measurable function. Then,  
$$m\mbox{-a.s.} \ \ x \ : \ \ \ \ (T(gf))(x)=(Tg)(x)f(x)$$
holds. 

{\rm (b)} \ $$m\mbox{-a.s.} \ \ x \ : \ \ \ \ (Tf)(x)=f(x)$$
holds. 

{\rm (c)} \ $$\int_{x\in X} ((1-T)f)(x)f(x)m(dx)=0$$ 
holds. 

{\rm (d)} \ $$\sigma_m\mbox{-a.s.} \ \ (x_0,x_1) \ : \ \ \ \ f(x_0)=f(x_1)$$
holds.} 

{\it Proof:} \ 
(1) \ Suppose that (a) holds. We show (b). 
For any $\mathbb R$-valued bounded $\cal B$-measurable function $g$, 
$\int_{x\in X} (Tf)(x)g(x) m(dx)
=\int_{x\in X} (T(1f))(x)$

\noindent 
$g(x) m(dx)
=\int_{x\in X} (T1)(x)f(x)g(x) m(dx)
=\int_{x\in X} f(x)g(x) m(dx)$ 
holds. 

(2) \ [(b)$\Rightarrow$(c)] is easy. 

(3) \ From Lemma 7, [(c)$\Rightarrow$(d)] holds. 

(4) \ Suppose that (d) holds. We show (a). 
Let $g$ be a $\mathbb R$-valued bounded $\cal B$-measurable function. 
Then, for any $\mathbb R$-valued bounded $\cal B$-measurable function $h$, 
$\int_{x\in X} h(x) (T(gf))(x) m(dx)
=\int_{(x_0,x_1)\in X^2} h(x_0)
g(x_1)f(x_1) \sigma_m(d(x_0,x_1))
=\int_{(x_0,x_1)\in X^2} h(x_0)g(x_1)f(x_0) \sigma_m(d(x_0,x_1))
=\int_{x\in X}h(x) (Tg)(x) f(x) m(dx)$ 
holds. 
\\ 
\hfill $\blacksquare$

{\bf Lemma 9:} \ 
{\it Let $m\in\cal G$. Let $\rho$ be a $[0,+\infty)$-valued bounded $\cal B$-measurable function. 
Then, $\rho m\in \cal G$ holds, if and only if $\rho m\in \cal I$ holds.}

{\it Proof:} \ 
(1) \ Suppose that $\rho m\in \cal G$ holds. 
We show that $\rho m\in \cal I$ holds. 
For any $\mathbb R$-valued bounded $\cal B$-measurable function $f$, 
$\int_{x\in X}(Tf)(x)\rho(x)m(dx)
=\int_{x\in X}(T1)(x)f(x)\rho(x)m(dx)=\int_{x\in X}f(x)\rho(x)m(dx)$ 
holds.

(2) \ Suppose that $\rho m\in \cal I$ holds. 
We show that $\rho m\in \cal G$ holds. 
For any $\mathbb R$-valued bounded $\cal B$-measurable function $f$, 
$\int_{x\in X}f(x)(T\rho)(x)m(dx)
=\int_{x\in X}(Tf)(x)\rho(x)m(dx)
=\int_{x\in X}f(x)\rho(x)m(dx)$ holds. 
[$m\mbox{-a.s.} \, x: \ (T\rho)(x)=\rho(x)$] holds. 
From Lemma 8, for any $\mathbb R$-valued bounded $\cal B$-measurable function $g$, 
[$m\mbox{-a.s.} \, x : \ (T(g\rho))(x)=(Tg)(x)\rho(x)$] 
holds. For any $\mathbb R$-valued bounded $\cal B$-measurable functions $f$ and $g$, 
$\int_{x\in X}(Tf)(x)g(x)\rho(x)m(dx)
=\int_{x\in X}f(x)(T(g\rho))(x)m(dx)
=\int_{x\in X}f(x)(Tg)(x)\rho(x)m(dx)$ 
holds. 
\hfill $\blacksquare$

{\bf Lemma 10:} \ 
{\it ${\cal G}_e\subset{\cal I}_e$ holds.}

{\it Proof:} \ 
Suppose that $m\in {\cal G}_e$ holds. We show that $m\in {\cal I}_e$ holds. 

(1) \ For any $\mathbb R$-valued bounded $\cal B$-measurable function $f$, 
$\int_{x\in X}(Tf)(x)m(dx)$

\noindent 
$=\int_{x\in X}(T1)(x)f(x)m(dx)
=\int_{x\in X}f(x)m(dx)$ 
holds. 

(2) \ Suppose that $t\in(0,1), \, m_0,m_1\in {\cal I}$ and $m=(1-t)m_0+tm_1$ hold. 
Then, $m_0\leq \frac{1}{1-t}m$ and $m_1\leq \frac{1}{t}m$ hold. 
There exist $[0,+\infty)$-valued bounded $\cal B$-measurable functions 
$\rho_0$ and $\rho_1$ such that $m_0=\rho_0m$ and $m_1=\rho_1m$ hold. 
From Lemma 9, $m_0, m_1\in\cal G$ holds. 
$m_0=m_1$ holds. 
\hfill $\blacksquare$

{\it Proof of Main Result:} \ 
${\cal R}={\cal G}$ holds. So, from Lemma 10, ${\cal R}_e={\cal G}_e\subset{\cal I}_e$ holds. 
\hfill $\blacksquare$

The following remark seems to be similar to the proposition (4.3.5) in [2]. 
Lemma 7 seems to be similar to the lemma (4.4.3) in [2]. 
As stochastic Ising model seems to be a typical example of symmetric Feller process, 
a reversible measure and a Gibbs measure are equivalent from the theorem (4.2.14) in [2].
While the corollary (4.4.20) in [2] asserts that an extreme Gibbs measure 
is an extreme invariant measure, 
it seems that (4.3.5) and (4.4.3) were essential for the proof in [2]. 

{\bf Remark:} \ 
{\it Let $m\in\cal G$. 
Let $\rho$ be a $[0,+\infty)$-valued bounded $\cal B$-measurable function. 
Suppose that $\int_{x\in X}\rho(x)m(dx)=1$ holds. 
Then, $\rho m\in \cal G$ holds, if and only if 
$[ \ \sigma_m\mbox{-a.s.} \, (x_0,x_1): \ \rho(x_0)=\rho(x_1) \ ]$ 
holds.} 

{\it Proof:} \ From Lemma 8, it is shown. 
(We do not go into details, as it is not difficult.) 
\hfill $\blacksquare$

\newpage 

\noindent
Reference:

\noindent
[1] \ Fukushima, M., Oshima, Y., Takeda, M., 
{\it Dirichlet forms and symmetric Markov processes}, 
Walter de Gruyter \& Co., Berlin, 1994. 

\noindent
[2] \ Liggett, T. M., {\it Interacting particle systems}, 
Springer-Verlag, New York, 1985.

\end{document}